\documentclass[12pt,reqno]{amsart}
\usepackage{amsfonts}
\usepackage{amssymb}
\usepackage{amsthm}
\usepackage{upref}
\usepackage{enumerate}

\usepackage{amscd}

\makeatletter 

\def\LaTeX{\leavevmode L\raise.42ex
    \hbox{\kern-.3em\size{\sf@size}{0pt}\selectfont A}\kern-.15em\TeX}
 \makeatother
 
\DeclareMathOperator{\clos}{clos}

\numberwithin{equation}{section}

\newtheorem{lemma}{Lemma}[section]
\newtheorem{theorem}[lemma]{Theorem} 
\newtheorem{corollary}[lemma]{Corollary}
\newtheorem{proposition}[lemma]{Proposition}

\theoremstyle{definition}

\newtheorem{example}[lemma]{Example}

\newtheorem{remark}[lemma]{Remark}

  \newcommand{\e}{\eqref}

\newcommand{\q}{\quad}

\newcommand{\ov}{\overline}
\newcommand{\wt}{\widetilde}

\newcommand{\ti}{\tilde}

\newcommand{\la}{\langle}
\newcommand{\ra}{\rangle}

\renewcommand\Im{\operatorname{Im}}

\newenvironment{pf}{\begin{proof}}{\end{proof}}

\def\qqq{\mathrel{\subset\mkern-15mu\lower.38ex\hbox{${\scriptscriptstyle\rightarrow}$}}}

\let\cal\mathcal

\let\Bbb\mathbb
 
\begin{document}
\title 
%{Self-adjoint Jacobi operators in the indeterminate case} 
%[ Self-adjoint Jacobi operators ]
%{Self-adjoint Jacobi operators associated with indeterminate moment problems}
{Self-adjoint Jacobi operators in the limit circle case}
\author{ D. R. Yafaev  }
\address{   Univ  Rennes, CNRS, IRMAR-UMR 6625, F-35000
    Rennes, France, SPGU, Univ. Nab. 7/9, Saint Petersburg, 199034 Russia, and     NTU Sirius, Olympiysky av. 1, Sochi, 354340 Russia}
\email{yafaev@univ-rennes1.fr}

% Asymptotic  behavior of orthogonal polynomials  with increasing recurrence coefficients in a critical case

% Asymptotic  behavior of orthogonal polynomials  with/for a critical growth  of recurrence coefficients
 

%\today{}
%   \date{9 April 2021}
\thanks {Supported by  project   Russian Science Foundation   17-11-01126}

\begin{abstract}
We consider symmetric Jacobi operators  with recurrence coefficients such that the corresponding difference equation is in the limit circle case. Equivalently, this means that the
 associated moment problem is indeterminate. Our main goal is to find a representation  for the resolvents of self-adjoint realizations $J$ of such Jacobi operators. This  representation implies  the classical Nevanlinna formula for the Cauchy-Stieltjes transforms of the spectral measures of the   operators $J$. We also efficiently describe domains of the operators $J$   in terms of boundary conditions at infinity.  
    \end{abstract}
    
    % In particular, we recover the classical Nevanlinna result stating that all operators $J$ have discrete spectra.

   %  in terms of asymptotic behavior of sequences $u_{n}$ in their domains as $n\to\infty$.
     
%\thispagestyle{empty}

 \maketitle

%\part{Introduction}
%***********************************************************

\section{Introduction}

\subsection{Basic definitions}

Given two sequences $a_{n}>0$ and $b_{n}=\bar{b}_{n}$ where $n\in {\Bbb Z}_{+}=\{0,1,\ldots \}$, 
 a 
 Jacobi matrix is defined by an equality
  \begin{equation}
{\cal J} = 
\begin{pmatrix}
 b_{0}&a_{0}& 0&0&0&\cdots \\
 a_{0}&b_{1}&a_{1}&0&0&\cdots \\
  0&a_{1}&b_{2}&a_{2}&0&\cdots \\
  0&0&a_{2}&b_{3}&a_{3}&\cdots \\
  \vdots&\vdots&\vdots&\ddots&\ddots&\ddots
\end{pmatrix} .
\label{eq:ZP+}\end{equation}
Let     $u= (u_{0}, u_{1}, \ldots)^\top =: (u_{n})$ be a column. 
  Then 
\begin{equation}
( {\cal J} u) _{0} =  b_{0} u_{0}+ a_{0} u_{1} \q \mbox{and}\q
( {\cal J} u) _{n} = a_{n-1} u_{n-1}+b_{n} u_{n}+ a_{n} u_{n+1} \q \mbox{for}\q n\geq 1.  
 \label{eq:ZP+1}\end{equation}

Let us consider a second-order difference equation  
\begin{equation}
  a_{n-1} u_{n-1}+b_{n} u_{n}+ a_{n} u_{n+1} =z u_{n}, \q n\geq 1,  
  \label{eq:Jy}\end{equation}
  associated with the operator $\cal J$. Clearly,   values  $u_{0}$ and $u_{1 }$  determine uniquely a solution   $u_{n}$ of  equation \e{eq:Jy}. In particular, 
the   solutions $p_{n}(z)$ and $q_{n}(z)$  are distinguished by  boundary conditions
 \begin{equation}
p_{0} (z)= 1,\q  q_{0} (z)= 0, 
  \label{eq:Pz}\end{equation}
  and
  \begin{equation}
p_1 (z)= (z-b_{0} )/a_{0} , \q  q_1 (z)= 1/a_{0} .
  \label{eq:Qz}\end{equation}
  It is easy to see that $p_{n}(z)$ is a  polynomial of degree $n$ and $q_{n}(z)$ is a  polynomial of degree $n-1 $.    Let us set $p(z)= (p_{0}(z),p_{1}(z),\ldots)^\top$, $q(z)= (q_{0}(z),q_{1}(z),\ldots)^\top$
  and    $e_{0}= (1,0,\ldots)^\top$.  Then  
  \begin{equation}
  {\cal J} p(z) =z p(z) \q \mbox{and} \q   {\cal J} q(z) =z q(z)  +e_{0}.
  \label{eq:PQzz}\end{equation}
  
   Let $u= ( u_{n} )$ and $v=(v_{n})$ be two solutions of equation \e{eq:Jy}. A direct calculation shows that their Wronskian
  \begin{equation}
\{ u, v  \} : = a_{n}  (u_{n}  v_{n+1}-u_{n+1}  v_{n})
\label{eq:Wr}\end{equation}
does not depend on $n\in{\Bbb Z}_{+}$. Clearly, the Wronskian $\{ u, v  \} =0$ if and only if the solutions $u$ and $v$ are proportional.  Note that $\{p,q \}=1$ so that
  \begin{equation}
 p_{n}(z) q_{n+1}(z)- p_{n+1}(z) q_{n}(z)=a_n^{-1}, \q \forall n\in {\Bbb Z}_{+}.
  \label{eq:PQW}\end{equation}

      We consider Jacobi operators $J$ acting    in 
 the space $\ell^2 ({\Bbb Z}_{+})$ with the scalar product   $\la u, v \ra=\sum_{n\in {\Bbb Z}_{+}}u_{n} \bar{v}_{n}$.
 Let $\cal D  $ be a dense in  $\ell^2 ({\Bbb Z}_{+})$ set of vectors $u $ with only a  finite number of non-zero components  $u_{n}$.
  The minimal Jacobi operator $J_{\rm min}$ is  defined   by the equality $J_{\rm min} u= {\cal J} u$ on domain 
 $  {\cal D} (J_{\rm min}) =\cal D$.
  It extends to a bounded operator defined on the whole space $\ell^2 ({\Bbb Z}_{+})$ 
  if and only if the corresponding recurrence coefficients $( a_{n} )\in \ell^\infty ({\Bbb Z}_{+}) $ and $ ( b_{n} )\in \ell^\infty ({\Bbb Z}_{+}) $.  In the general case, we introduce also
 the maximal  operator $J_{\rm max}$   by the   same  formula $J_{\rm max} u= {\cal J} u$ on a set $ {\cal D} (J_{\rm max})$  of all vectors $u \in \ell^2 ({\Bbb Z}_{+})$ such that ${\cal J}u \in \ell^2 ({\Bbb Z}_{+})$. 
    The operator $J_{\rm min} $ is
symmetric  in the space $\ell^2 ({\Bbb Z}_{+})$, and its adjoint operator  $J^*_{\rm min} = J_{\rm max}$.

The operator $J_{\rm max}$  is not of course symmetric.  
     For all $u, v \in {\cal D} (J_{\rm max})$, we have an identity  (the Green formula)
       \begin{equation}
\la {\cal J} u,v \ra - \la u,  {\cal J} v \ra = \lim_{n\to\infty}  a_{n } (u_{n +1}\bar{v}_{n }- u_{n }\bar{v}_{n +1})
\label{eq:qres4}\end{equation}
where the limit in the right-hand side exists.
 Indeed, a direct calculation shows that
      \[
\sum_{m=0}^n ( {\cal J} u)_{m}\bar{v}_{m} -\sum_{m=0}^n u_{m} ( {\cal J} \bar{v})_{m} =     a_{n } (u_{n +1}\bar{v}_{n }- u_{n }\bar{v}_{n +1}).
\]
Passing here to the limit $n\to\infty$ and using that ${\cal J} u  \in \ell^2 ({\Bbb Z}_{+})$ for $u \in {\cal D} (J_{\rm max})$, we obtain \e{eq:qres4}.

 The comprehensive presentation of the theory of self-adjoint Jacobi operators can be found in  the books \cite{AKH} (Chapters~I and II), \cite{Ber} (Chapter~VII), \cite{Schm} (Part~I)  and the survey \cite{Simon}.

  \subsection{Determinacy versus indeterminacy}
  
    Since the operator $J_{\rm min}  $  commutes with the complex conjugation,  its deficiency indices 
   \[
   d_{\pm}: =\dim\ker (J_{\rm max}-z I),\q \pm \Im z>0,
   \]
     are equal, i.e. $d_{+} =d_{-}=:d $, and, so, $J_{\rm min}  $ admits self-adjoint extensions.  Here and below $I$ is the identity operator. For an arbitrary $z\in{\Bbb C}$,
       let us consider the equation ${\cal J}u=zu$. By definition \e{eq:ZP+1}, its solutions  are given by the formula $u=u_{0}p(z)$. They belong to $ {\cal D} (J_{\rm max})$ if and only if $p (z) \in  \ell^2 ({\Bbb Z}_{+})$.  Therefore $d=0$ if $p (z)\not\in  \ell^2 ({\Bbb Z}_{+})$  for $\Im z\neq 0$; otherwise $d=1$.  
     
      This dichotomy can also be stated in terms of the Weyl limit point/circle theory and in terms of the determinacy/indeterminacy  of the associated moment problem.  Recall that,
   similarly to differential equations,   difference equation \e{eq:Jy} always has a non-trivial solution in $\ell^2 ({\Bbb Z}_{+})$ for $\Im z \neq 0$. This solution is either unique (up to a constant factor) or all solutions of \e{eq:Jy} belong to $\ell^2 ({\Bbb Z}_{+})$.
The first instance is known as the limit point (LP) case and the second one --  as the limit circle (LC) case.   In the LP case $p (z)\not\in  \ell^2 ({\Bbb Z}_{+})$   and in the LC case $p (z) \in  \ell^2 ({\Bbb Z}_{+})$  for $\Im z\neq 0$. Thus, the operator $J_{\rm min}$ is essentially self-adjoint if and only if the LP case occurs.
 
 % Clearly, equation \e{eq:Jy} complemented by the boundary condition $  b_{0}u_{0}+ a_{0}u_{1} =zu_{0}$ has a unique (up to a constant factor) solution $ u_{n}=P_{n}  (z)$. It satisfies the equation $J_{\max}u=z u$ if and only if $\{P_{n}  (z)\} \in \ell^2 ({\Bbb Z}_{+})$. It follows that the deficiency indices of the operator $J_{\min}$  are either $(0,0)$ or $(1,1)$.

%     In the framework of the Weyl limit point/circle theory, 
% the first instance is known as the limit point (LP) case and the second one --  as the limit circle (LC) case. 
% In the LP case the operator $J_{\rm min}$ is essentially self-adjoint  while it has an one-parameter family of self adjoint extensions  in the LC case.

  %   This dichotomy can also be stated in terms of the limit point/circle theory.  Recall that,
 %  similarly to differential equations, the  difference equation \e{eq:Jy} always has a non-trivial solution in $\ell^2 ({\Bbb Z}_{+})$ for $\Im z \neq 0$. This solution is either unique (up to a constant factor) or all solutions of \e{eq:Jy} belong to $\ell^2 ({\Bbb Z}_{+})$.

   %    The operator $J_{\rm min}  $ may be  essentially  self-adjoint, that is,  its closure $\clos J_{\rm min}=J_{\rm max}$,   ({\it determinate}  case) or have multiple self-adjoint extensions ({\it indeterminate}  case).
   % This dichotomy can also be stated in terms of solutions of the difference equations \e{eq:Jy} and in terms of the determinacy of the associated moment problems.

As is well known (see, for example, Theorem~6.16 in \cite{Schm}),  in the LC case
inclusions
  \begin{equation}
  p(z)\in  \ell^2 ({\Bbb Z}_{+}) \q \mbox{and} \q    q (z)\in  \ell^2 ({\Bbb Z}_{+})  
  \label{eq:PQz}\end{equation}
hold true for all $z\in {\Bbb C}$.  It now follows from \e{eq:PQW} that  in the LC case  necessarily  
    \begin{equation}
\sum_{n=0}^\infty a_{n}^{-1}<\infty.
\label{eq:noCarl}\end{equation}
Equivalently, this fact can be stated   
  as a sufficient condition (known as the Carleman condition) 
     \begin{equation}
\sum_{n=0}^\infty a_{n}^{-1}=\infty
\label{eq:Carl}\end{equation}
for the essential self-adjointness of the operator $J_{\rm min}$. 

In general,  the essential self-adjointness of  operators $J_{\rm min}$ is determined by both sets of coefficients $(a_{n})$ and $(b_{n})$. This is briefly discussed in Sect.~3. Here we only note that, even in the case $(b_{n})=0$,
condition \e{eq:Carl} is not necessary for the essential self-adjointness of   $J_{\rm min}$. However Theorem~1.5 in Chaper~VII of the book \cite{Ber} shows that this is   true under an additional assumption 
  $a_{n-1}a_{n+1}\leq a_{n}^2$
  %\label{eq:ber}\end{equation}
(for large $n$). Without   this assumption, the operators  $J_{\rm min}$ may be  essentially self-adjoint even for the coefficients $(a_{n})$ satisfying  conditions \e{eq:noCarl}  and $(b_{n})=0$.
 (see \cite{Kost}).
 
 We also note a link with the  Hamburger moment problem
      \begin{equation}
s_{n} = \int_{-\infty}^\infty \lambda^n d\rho(\lambda).
\label{eq:moment}\end{equation}
Recall that these equations for a nonnegative  measure $d\rho(\lambda)$ with infinite support have a solution if and only if the sequence $(s_{n})$ is positive definite, i.e., 
    \begin{equation}
\sum_{n,m\in{\Bbb Z}_{+}}s_{n+m}\xi_{m}\bar{\xi_{n}}>0\q \forall \xi=(\xi_{n})\in {\cal D}, \;\xi\neq 0.
\label{eq:moment1}\end{equation}
This problem is called determinate if the measure satisfying 
\e{eq:moment} is unique. Otherwise it  is called indeterminate.

For a Jacobi matrix \e{eq:ZP+}, we  set
$    s_{n}=({\cal J}^n e_{0}, e_{0})$.
 Then condition \e{eq:moment1} is satisfied. If  $J$ is an arbitrary self-adjoint extension of the operator $J_{\rm min}$  and $E_{J}(\lambda)$ is its spectral family,  then the     measure $d\rho_{J}(\lambda)= d(E_{J}(\lambda)e_{0}, e_{0})$ satisfies equation \e{eq:moment}.
  It turns out (see Theorem~2 in \cite{Simon})  that  moment problem \e{eq:moment} is   determinate  if and only if the operator $J_{\rm min}$ is     essentially  self-adjoint.

We also note that the  essential  self-adjointness of an operator $J_{\rm min}$ is determined by a behavior of its coefficients $a_{n}$ and $b_{n}$ for large $n$ only. Indeed, if  $( \ti{a}_{n}-a_{n})\in \ell^\infty ({\Bbb Z}_{+}) $  and $( \ti{b}_{n}- b_{n})\in \ell^\infty ({\Bbb Z}_{+}) $, then the operator $\wt{J}_{\rm min}-J_{\rm min}$ is bounded so that the deficiency indices of the operators    $\wt{J}_{\rm min}$ and $J_{\rm min}$ are the same.

% It follows that Jacobi operators $J$ and $\wt{J}$ are  both in the determinate (or indeterminate) case if the corresponding Jacobi coefficients satisfy the conditions $\{a_{n}-\ti{a}_{n}\}\in \ell^\infty ({\Bbb Z}_{+}) $  and $\{b_{n}-\ti{b}_{n}\}\in \ell^\infty ({\Bbb Z}_{+}) $.

    \subsection{Self-adjoint extensions}
    
  Recall that
\[
\clos J_{\min}= J_{\min}^{**}  = J_{\max}^{*}. 
\]
  In this paper we are interested
      in the LC case where the minimal Jacobi operators  $J_{\rm min}$  are not essentially self-adjoint so that      \begin{equation}
\clos J_{\rm min}\neq J_{\rm max}=J_{\rm min}^*.
  \label{eq:cl}\end{equation}

 Let us set
   \begin{equation}
{\sf D }:={\cal D}(\clos J_{\rm min}).
  \label{eq:MIN}\end{equation}
  For a vector $h\in\ell^2 ({\Bbb Z}_{+})$, we denote by $\{h\}$ the one dimensional subspace of $\ell^2 ({\Bbb Z}_{+})$ spanned by the vector $h$. The symbol $\dotplus$ denotes the direct sum of subspaces. The following description  (see Lemma~6.22 and Theorem~6.23 in \cite{Schm}  or Theorem~2.6 in  \cite{Simon}) of self-adjoint extensions of the operators $J_{\rm min}$ is a modification of von Neumann formulas adapted to Jacobi operators. 
  
   \begin{proposition}\label{Neum}
    Let relation \e{eq:cl} be satisfied. Choose some $\zeta\in {\Bbb R}$.
Then 
\begin{equation}
{\cal D}(J_{\rm max})= {\sf D }\dotplus \{p(\zeta)\}\dotplus \{q(\zeta)\}.
  \label{eq:Neum}\end{equation}
  All self-adjoint extensions $J_{t}$ of the operator $J_{\rm min}$ are parametrized by numbers $t\in{\Bbb R}$ and
  $t=\infty$. The domains of the operators $J_{t}$ are determined by the equalities
\begin{equation}
{\cal D}(J_t)= {\sf D }\dotplus \{t p(\zeta) +q(\zeta)\}, \q t\in{\Bbb R},
    \label{eq:Neum1}\end{equation}
  and
\begin{equation}
{\cal D}(J_\infty)= {\sf D }\dotplus \{ p(\zeta)  \}.
  \label{eq:Neum2}\end{equation}
 \end{proposition} 
 
 Below we use Proposition~\ref{Neum}  for $\zeta=0$ only.
   
 Our main goal is to find an explicit formula for the resolvents $R_{t} (z)= (J_{t}-zI)^{-1}$  of   self-adjoint operators $J_{t}$.  To that end, we define, for all $z\in{\Bbb C}$, a bounded operator ${\cal R} (z)$ in the  space $\ell^2 ({\Bbb Z}_{+})$  by an equality
        \begin{equation}
 ( {\cal R} (z)h)_{n} =  q_{n}(z) \sum_{m=0}^n  p_{m} (z) h_{m}+   p_{n}(z)   \sum_{m=n+1}^\infty  q_{m} (z) h_{m}.\label{eq:RR11}\end{equation} 
 Actually, the operator ${\cal R} (z)$ belongs to the Hilbert-Schmidt class.
We prove (see  Theorem~\ref{res}) that, in  a natural sense,  ${\cal R} (z)$ can be considered as a quasiresolvent of the operator $J_{\rm max}$. A representation for the resolvents $R_{t} (z)$ is obtained in Theorem~\ref{RES}.
 
 Another goal of the paper is to find an alternative description of self-adjoint extensions  of the operator $J_{\rm min}$ in terms of an asymptotic behavior of solutions $u_{n}$ of the Jacobi equation  \e{eq:Jy} for $n\to\infty$.  This is discussed in Section~3.

\section{Resolvents of self-adjoint Jacobi operators}

 In this section we only suppose that a minimal Jacobi operator $J_{\rm min}$ is not essentially self-adjoint so that we are in the LC case.

      \subsection{Quasiresolvent}

Recall that in the LC  case
 inclusions \e{eq:PQz}  are satisfied for all $z\in {\Bbb C}$.    Let us define an operator ${\cal R}  (z)$ playing the role of the resolvent of the operator $J_{\rm max}$ by   equality \e{eq:RR11}.  In particular, it follows from \e{eq:Pz}, \e{eq:Qz} that  
  \begin{equation}
  ({\cal R} (z)h)_{0}=\la h, q(\bar{z})\ra 
  \label{eq:r01}\end{equation} 
and
 \begin{equation}
   ({\cal R} (z)h)_{1}= h_{0}  a_{0}^{-1}+\la h, q(\bar{z})\ra  (z-b_{0})a_{0}^{-1} 
\label{eq:r02}\end{equation} 
for all  $h\in {\ell}^2  ({\Bbb Z}_{+})$. 
 In view of inclusions \e{eq:PQz} the operators ${\cal R}  (z)$ are bounded in the space $\ell^2 ({\Bbb Z}_{+})$ for all $z\in {\Bbb C}$.  Note also that ${\cal R} (z)^*={\cal R} (\bar{z})$.

A proof of the following statement is  close to the construction (see, e.g., Lemma~5.1 in \cite{JLR}) of the resolvents for  essentially self-adjoint Jacobi operators.    
   
    \begin{theorem}\label{res}
    Let relation \e{eq:cl} be satisfied.
 For all $z \in {\Bbb C}$, we have
  \begin{equation}
{\cal R}  (z): \ell^2 ({\Bbb Z}_{+})\to  {\cal D} (J_{\rm max} )
\label{eq:qres}\end{equation}
and
  \begin{equation}
(J_{\rm max} -zI) {\cal R}  (z)=I.
\label{eq:qres1}\end{equation}
   \end{theorem} 
 
  \begin{pf} 
  Recall that an operator ${\cal J}$  was defined by equalities \e{eq:ZP+1}. We will  check that 
    \begin{equation}
  (  ({\cal J}-z I) {\cal R} (z)h)_{n}=h_{n} 
    \label{eq:RRm}\end{equation}
for all $n \in {\Bbb Z}_{+}$ and $h=\{h_{n}\}\in \ell^2 ({\Bbb Z}_{+})$.  For $n=0$, we have
\[
(  ({\cal J}-z I) {\cal R} (z)h)_0= (b_{0}-z)({\cal R} (z)h)_0  + a_{0}({\cal R} (z)h)_1=h_{0}
\]
according to formulas \e{eq:r01} and \e{eq:r02}. For $n\geq 1$, we rewrite definition  \e{eq:RR11}   as
 \[
 ( {\cal R} (z)h)_{n} =  q_{n} (z) x_{n}  (z) +   p_{n} (z)   y_{n} (z)
 \]
where
   \begin{equation}
x_{n}(z) =\sum_{m=0}^n  p_{m} (z) h_{m},\q    y_{n}(z) =\sum_{m=n+1}^\infty  q_{m} (z) h_{m}.
\label{eq:RR2}\end{equation} 
  It now follows from definition \e{eq:ZP+1}     that
   \begin{multline}
 ( ({\cal J} -z  I) {\cal R} (z) h)_{n} =  a_{n-1}\big(q_{n-1}x_{n-1}+p_{n-1} y_{n-1} \big)
  \\
  + (b_{n}-z) \big(q_{n}x_{n}+ p_{n} y_{n}  \big)+ a_{n}\big(q_{n+1}x_{n+1} +p_{n+1} y_{n+1}    \big), \q n\geq 1.
\label{eq:RR3}
\end{multline}
 According to \e{eq:RR2}  we have
\[
q_{n-1}x_{n-1}+p_{n-1} y_{n-1}=q_{n-1}(x_{n}-p_{n} h_{n})+  p_{n-1}(y_{n}+q_{n} h_{n}) 
\]
and
 \begin{multline*}
q_{n+1} x_{n+1} +p_{n+1} y_{n+1} =
q_{n+1}(x_{n}  +p_{n+1}h_{n+1}) 
\\
 +  p_{n+1} (y_{n} - q_{n+1}h_{n+1})=
q_{n+1}x_{n}  +  p_{n+1}y_{n} .
\end{multline*}
 Substituting these expressions into the right-hand side of \e{eq:RR3},  we see that
 \begin{multline*}
  \big( ({\cal J} -z I) {\cal R} (z) h\big)_{n} =  a_{n-1}\Big(q_{n-1}(x_{n}- p_{n} h_{n})+  p_{n-1}(y_{n}+q_{n} h_{n}) \Big)
  \\
  + (b_{n}-z) \big(q_{n}x_{n}+ p_{n} y_{n}  \big)+ a_{n}\big(q_{n+1}x_{n}  +  p_{n+1}y_{n}    \big).
\end{multline*}
Let  us now collect together all terms containing $x_{n}$, $y_{n}$ and $h_{n}$.  Then
 \begin{multline}
  \big( ({\cal J} -z I) {\cal R} (z) h\big)_{n} = \Big(a_{n-1}q_{n-1}+  (b_{n}-z) q_{n} + a_{n} q_{n+1}\Big) x_{n}
  \\
    + \Big(a_{n-1}p_{n-1}+  (b_{n}-z)p_{n} + a_{n}p_{n+1}\Big) y_{n}+a_{n-1}\big( -p_{n} q_{n-1}+  p_{n-1} q_{n}\big) h_{n}  .
\label{eq:RR3A}
\end{multline}
The coefficients at $x_{n}$ and $y_{n}$ equal zero by virtue of   Jacobi equation \e{eq:Jy}  for $( q_{n} )$ and $(p_{n})$, respectively. Since $\{p,q \} =1$, the right-hand side of 
\e{eq:RR3A} equals $ h_{n}$.   
This proves \e{eq:RRm}  whence $   ({\cal J}-z I) {\cal R} (z)h=h $. In particular, we see that ${\cal R}  (z) h\in {\cal D} (J_{\rm max})  $ so that ${\cal J}$ can be replaced here by $J_{\rm max}$.
This yields both \e{eq:qres} and \e{eq:qres1}. 
    \end{pf}
    
    \begin{remark}\label{res3}
    In definition  \e{eq:RR11},
    one can replace   $q_{n}(z)$ by the polynomials $\wt{q}_{n}(z)= q_{n}(z)+ c p_{n}(z)$ for an arbitrary  $c\in {\Bbb C}$. Then $\wt{\cal R} (z)= {\cal R} (z) +c \la \cdot, p(\bar{z})\ra p(z)$ and formulas \e{eq:qres}, \e{eq:qres1}  remain true.
      \end{remark}

    Since $p(z)$ is a unique (up to a constant factor)  solution of the homogeneous equation $(J_{\rm max} -zI) u =0$, we can also state 
    
     \begin{corollary}\label{res1}
     All solutions of the equation
       \[
(J_{\rm max} -zI) u =h \q \mbox{where}  \q z\in {\Bbb C} \q \mbox{and}\q h\in \ell^2 ({\Bbb Z}_{+})
\]
  for $u\in {\cal D} (J_{\rm max} )$ are  given by the formula
 \begin{equation}
  u = \Gamma p(z)+ {\cal R}  (z) h \q \mbox{for some} \q \Gamma=\Gamma (z;h) \in{\Bbb C}.
\label{eq:qres3}\end{equation}
   \end{corollary} 
   
 %  Alternatively, this result can be stated in the following way.
   
 %  \begin{corollary}\label{res1X}
  % Choose some $z\in {\Bbb C}$.  All vectors $u\in {\cal D} (J_{\rm max} )$ admit representations
  %  \begin{equation}
%  u =  \boldsymbol{\Gamma}  (z; u) p(z)+ {\cal R}  (z) (J_{\rm max} -z I) u
%\label{eq:qres3X}\end{equation}
%with a complex number  $ \boldsymbol{\Gamma}  (z; u)$ depending on $z$ and $u$.
 %  \end{corollary} 

An  asymptotic relation for $ ( {\cal R} (z)h)_{n}$ is a direct consequence of definition  \e{eq:RR11}
   and condition  \e{eq:PQz}:
      \begin{equation}
 ( {\cal R} (z)h)_{n} =  q_{n}(z) \la    h, p (\bar{z}) \ra+   o(|p_{n} (z)| +|q_{n}(z)| )\q {\rm as}\q n\to\infty.
 \label{eq:Ras}\end{equation}
 This relation   can be supplemented by the following result.
 
     \begin{proposition}\label{AS}
    For all $z\in {\Bbb C}$ and  all $h\in \ell^2 ({\Bbb Z}_{+})$, we have
    \begin{equation}
 u:= {\cal R} (z)h - q (z) \la    h, p (\bar{z}) \ra\in {\sf D}.
 \label{eq:Ras1}\end{equation}
        \end{proposition} 
        
              \begin{pf}
              Let first  $h\in {\cal D}$. Then $ ( {\cal R} (z)h)_{n} =  q_{n}(z) \la    h, p (\bar{z}) \ra$ for sufficiently large $n$ so that $u\in {\cal D}\subset {\sf D}$.

Let now $h$ be an arbitrary vector in  $  \ell^2 ({\Bbb Z}_{+})$.     Observe that $ u \in {\sf D}$ if and only if there exists a sequence $u^{(k)}\in {\cal D}$ such that 
                \begin{equation}
u^{(k)}\to u \q \mbox{and} \q {\cal J}u^{(k)}\to {\cal J}u \q \mbox{as} \q k\to\infty.
 \label{eq:Ras2}\end{equation}
            Let us take any sequence $h^{(k)}\in {\cal D}$ such that   $h^{(k)}\to h$ and set
                 \[
 u^{(k)}= {\cal R} (z)h^{(k)} - q (z) \la    h^{(k)}, p (\bar{z}) \ra .
 \]
 Then $u^{(k)} \in {\cal D}$ and  $u^{(k)}\to u$ as $k\to\infty$ because the operator $ {\cal R} (z)$ is bounded.  It follows from equalities  \e{eq:PQzz}  and \e{eq:qres1}  that
 \[
({\cal J}  -z) u^{(k)}=  h^{(k)} - e_{0} \la    h^{(k)}, p (\bar{z}) \ra \to h - e_{0} \la    h , p (\bar{z} ) \ra =({\cal J}  -z) u 
\]
as $k\to\infty$. This proves relations \e{eq:Ras2} whence $u\in {\sf D}$. 
                              \end{pf}
                              
                              In view of Proposition~\ref{AS}  representation \e{eq:qres3}  is consistent with formula  \e{eq:Neum}.
      
          \subsection{Main result}
          
          Let us first find a link between the polynomials $p(z)$, $q(z)$ for an arbitrary  $z \in {\Bbb C}$ and  the polynomials $p(0)$, $q(0)$.  We always suppose that
               condition \e{eq:cl} is  satisfied which implies inclusions \e{eq:PQz} for all $z \in {\Bbb C}$.
          
          \begin{lemma}\label{PQ}
 For all $z \in {\Bbb C}$, we have
  \begin{equation}
p(z) =\big( 1-z \la p(z), q(0)\ra \big) p(0) + z {\cal R}(0) p(z)
\label{eq:PP}\end{equation}
and
  \begin{equation}
q(z) =    - z \la q(z), q(0)\ra  p(0) + q(0)+ z {\cal R}(0) q(z).
\label{eq:QQ}\end{equation}
 \end{lemma}
 
   \begin{pf}
   To prove \e{eq:PP}, we set
    \begin{equation}
   u=p(z) -z {\cal R}(0) p(z)
\label{eq:PP1}\end{equation}
so that 
    \begin{equation}
u_{0}= 1-z \la p(z), q(0)\ra
\label{eq:PP2}\end{equation} 
  according to  equality   \e{eq:r01} (for $z=0$).
   It follows from relation \e{eq:qres1}  that ${\cal J} u={\cal J} p(z) -z p(z) =0$ whence $u=u_{0} p(0)$. 
   Now definition \e{eq:PP1} implies that $ p(z)= u_{0} p(0)+ z {\cal R}(0) p(z)$. It remains to use \e{eq:PP2}.

   The proof of \e{eq:QQ}  is quite similar.  We now set
    \begin{equation}
   v=q(z) - q(0) -z {\cal R}(0) q(z)
\label{eq:QQ1}\end{equation}
so that  $v_{0}= -z \la q(z), q(0)\ra$.     It again follows from relations  \e{eq:PQzz} and  \e{eq:qres1} 
that  ${\cal J} v= {\cal J} q(z)- {\cal J} q(0)- z q(z)=0$  whence $v=v_{0} p(0)= - z \la q(z), q(0)\ra  p(0)$. Therefore 
\e{eq:QQ} is a direct consequence of definition \e{eq:QQ1}. 
\end{pf}
 
 Putting together Lemma~\ref{PQ}  with Proposition~\ref{AS} (for $z=0$), we can also state the following result.
 
   \begin{lemma}\label{PQ+} 
 For all $z \in {\Bbb C}$, we have
  \begin{equation}
p(z) -\big( 1-z \la p(z), q(0)\ra \big) p(0) - z \la p(z), p(0)\ra  q(0) \in {\sf D}
\label{eq:PP+}\end{equation}
and
  \begin{equation}
q(z) + z \la q(z), q(0)\ra  p(0) - \big( 1+ z \la q(z), p(0)\ra \big)   q(0)     \in {\sf D}.
\label{eq:QQ+}\end{equation}
 \end{lemma}

          Now it easy to construct     resolvents  of   self-adjoint extensions $J_{t}$ of the operators $J_{\rm min}$.  Recall that domains of the operators $J_{t}$ are defined by relations \e{eq:Neum1}  and \e{eq:Neum2}.
     
     \begin{theorem}\label{RES}
          Let assumption \e{eq:cl} hold.
     For all $z\in {\Bbb C}$ with $\Im z\neq 0$ and all $h\in\ell^2 ({\Bbb Z}_{+})$, the resolvent $R_t (z)= (J_t-zI)^{-1}$ of the operator $J_t$ is given by an  equality 
           \begin{equation}
R_{t} (z) h = \gamma_{t} (z) \la h,  p(\bar{z})\ra p(z)+ {\cal R}(z)h
      \label{eq:RES1}\end{equation}
      where  
           \begin{equation}
\gamma_{t}(z)= \frac{z \la q(z), q(0)\ra+\big( 1+ z \la q(z), p(0)\ra \big)t } 
{1- z \la p(z), q(0)\ra -z \la p(z), p(0)\ra t }, \q t\in {\Bbb R},
    \label{eq:RES2}\end{equation}
    and
            \begin{equation}
\gamma_{\infty}(z)= - \frac{ 1+ z \la q(z), p(0)\ra } 
{z \la p(z), p(0)\ra   }.
    \label{eq:RES3}\end{equation}
    \end{theorem} 
    
    \begin{pf}
    According to Theorem~\ref{res}  and Corollary ~\ref{res1}  a vector $u=R_{t}h$ is given by equality \e{eq:qres3}  where 
  $ \Gamma= \Gamma_{t}(z;h)$ is a bounded linear functional of $h\in \ell^2 ({\Bbb Z}_{+})$ so that
  $\Gamma_{t}(z;h)= \la h, \ti{p_{t}}(z)\ra$
   for some vector $ \ti{p}_{t} (z)\in \ell^2 ({\Bbb Z}_{+})$.
   Since $R_{t} (z)^*=R_{t} (\bar{z})$ and ${\cal R} (z)^*={\cal R}  (\bar{z})$,  we see that 
   \[
   \la h, \ti{p_{t}}(\bar{z})\ra p(\bar{z}) =\la h,  p (z)\ra \ti{p_{t}}(z)
   \]
   for all $h\in\ell^2 ({\Bbb Z}_{+})$. It follows that $\ti{p_{t}}(z)= \ov{\gamma_{t} (z)}p(\bar{z})$ for some $\gamma_{t} (z) \in{\Bbb C}$. This yields representation \e{eq:RES1} where the constant 
$\gamma_{t} (z) $ is determined by the condition $R_{t} (z) h\in {\cal D} (J_{t})$.
In view of relations  \e{eq:RES1}  and  \e{eq:Ras1}, this is equivalent to the condition
      \begin{equation}
R_{t} (z) h -\la h,  p(\bar{z})\ra \big( \gamma_{t} (z)  p(z) + q(z)\big)\in{\sf D}.
      \label{eq:RES4}\end{equation}
     
   Proposition~\ref{Neum} where $\zeta=0$ means that
       \begin{equation}
R_{t} (z) h - X  \big( t  p(0) + q(0)\big)\in{\sf D} \q \mbox{if} \q t\in {\Bbb R}
 \q \mbox{and} \q R_\infty (z) h - X  p(0)  \in{\sf D}
      \label{eq:REA}\end{equation}
      for some number $X=X_{t}(z)\in {\Bbb C}$.  Comparing  \e{eq:RES4} and \e{eq:REA}, we see that
      \begin{equation}
\gamma_{t} (z)  p(z) + q(z) - Y  \big( t  p(0) + q(0)\big)\in{\sf D} \q \mbox{if} \q t\in {\Bbb D}
      \label{eq:REA1}\end{equation}
 and
       \begin{equation}
 \gamma_\infty (z)  p(z) + q(z) - Y  p(0)  \in{\sf D}
      \label{eq:REA2}\end{equation}
         for some number $Y\in {\Bbb C}$. 
      
     Using     Lemma~\ref{PQ+}, we see that 
       inclusion \e{eq:REA1} is  equivalent to an equality
\begin{multline*}
  \gamma_{t} (z) (\big( 1-z \la p(z), q(0)\ra \big) p(0)  +z \la p(z), p(0)\ra q(0)\Big)
\\
 + \Big(- z \la q(z), q(0)\ra p(0) + \big( 1+ z \la q(z), p(0)\ra \big)   q(0)   \Big) = Y (t p(0) + q(0)) .
 \end{multline*}
  Comparing here the coefficients at $p(0)$ and $q(0)$, we obtain an equation
      \[
      \frac{  \gamma_{t} (z)  ( 1-z \la p(z), q(0)\ra )- z \la q(z), q(0)\ra }
      {\gamma_{t} (z) z \la p(z), p(0)\ra + 1+ z \la q(z), p(0)\ra } =t.
      \]
      Solving this equation with respect to $\gamma_{t} (z) $, we arrive at formula  \e{eq:RES2}.
         Similarly, substituting expressions \e{eq:PP+} and \e{eq:QQ+} into the left-hand side of \e{eq:REA2}, we see that inclusion \e{eq:REA2} holds true  if and only if the coefficient at $q(0)$   equals zero. This yields formula  \e{eq:RES3}.
        \end {pf}

          We emphasize that, for different $t$, the resolvents  $R_{t}(z)$ of the operators $J_{t}$ differ from each other only by the coefficient $\gamma_{t}(z)$ at the rank one operator $\la\cdot, p(\bar{z})\ra p(z)$.
   Observe also that
    $  \ov{\gamma_{t}(z)}= \gamma_t (\bar{z})$.

     Since  $\la {\cal R}(z) e_{0}, e_{0}\ra=0$, we see that $\la R_{t}(z) e_{0}, e_{0}\ra=\gamma_{t} (z)$. Thus Theorem~\ref{RES}  implies the classical Nevanlinna representation obtained in \cite{Nevan}  for the Cauchy-Stieltjes transform of the spectral measures $d\rho_{t} (\lambda)= d (E_{t} (\lambda) e_{0}, e_{0})$ of the operators $J_{t}$.
     
     \begin{corollary}\label{RESc}
     For all $z\in {\Bbb C}$ with $\Im z\neq 0$,  we have
     \[
     \int_{-\infty}^\infty
 (\lambda-z)^{-1} d\rho_{t}(\lambda)= \gamma_{t} (z),
 \]
 where $\gamma_{t} (z)$ is given by equalities \e{eq:RES2}  or  \e{eq:RES3}.
    \end{corollary} 

Since the functions $ \la p(z), p(0)\ra  $ and $ \la p(z), q(0)\ra  $  are entire, it follows from \e{eq:RES2}  and  \e{eq:RES3} that the spectra of the operators $J_{t}$  and $J_{\infty}$ are discrete. The eigenvalues of the operators $J_{t}$ are given by the equation 
   $z \la p(z), p(0)\ra t + z \la p(z), q(0)\ra =1$ if $t\in{\Bbb R}$ and  by the equation 
   $z \la p(z), p(0)\ra =0$ if $t=\infty$. This result is also due to R.~Nevanlinna.
   
     \begin{corollary}\label{RESb}
     If $z\in {\Bbb C}$  is  not an eigenvalue of the operator $J_{t}$, then its resolvent $R_{t}  (z)$ is in the Hilbert-Schmidt class.
    \end{corollary} 
    
     \begin{remark}\label{ze}
     Theorem~\ref{RES} remains true if one  uses an arbitrary real point $\zeta$ in parametrization of self-adjoint extensions of the operator $J_{\rm min}$ in Proposition~\ref{Neum}. In this case $p(0)$ and $q(0)$ in formulas 
     \e{eq:RES2}  and  \e{eq:RES3}  should be replaced by $p(\zeta)$ and $q(\zeta)$; the factor $z$ should be replaced by $z-\zeta$. 
       \end{remark}  
       
       Finally, let us compare resolvent formulas in the LP and LC cases. In the LP case the resolvent $R(z)$ of the Jacobi operator $J=\clos J_{\rm min}$   is given (see, e.g., Lemma~5.1 in \cite{JLR}) by the relation
        \begin{equation}
 ( R (z)h)_{n} = \frac{1}{\{p(z), g(z)\}} \Big(g_{n}(z) \sum_{m=0}^n  p_{m} (z) h_{m}+   p_{n}(z)   \sum_{m=n+1}^\infty  g_{m} (z) h_{m}\Big)
 \label{eq:R-LP}\end{equation}
 where $g(z)$ is a unique (up to a constant factor) solution of equation \e{eq:Jy}  belonging to $\ell^2 ({\Bbb Z}_{+})$. It can be chosen in a form $g(z)= q(z)+ w(z) p(z)$ where $w(z)$ is known as the Weyl function. Substituting this expression into \e{eq:R-LP}, we see that {\it formally}
   \begin{equation}
  R (z)  = w(z)\la\cdot,p(\bar{z}\ra p(z) + {\cal R} (z)
 \label{eq:R-LP1}\end{equation}
 where $ {\cal R} (z)$ is given by equality  \e{eq:RR11}. This relation looks algebraically similar to \e{eq:RES1}  where $\gamma_{t}(z)$ plays the role of the Weyl function $w(z)$. Note, however, that in the LP case 
 $w(z)$ is determined uniquely by the condition $q(z)+ w(z) p(z)\in \ell^2 ({\Bbb Z}_{+})$ while in the LC case 
 $\gamma_{t}(z)$ depends on the choice of a self-adjoint extension of the operator $  J_{\rm min}$.  We also emphasize that  relation \e{eq:R-LP1} is only formal  because $p(z)$ and $q(z)$ are not in $\ell^2 ({\Bbb Z}_{+})$.

  \section{Asymptotic behavior of solutions of the Jacobi equation}
  
  In this section we consider a specific class of Jacobi operators $J_{\rm min}$ that are not essentially   essentially self-adjoint.  We   describe  all their self-adjoint extensions $J$  in terms of asymptotic behavior of solutions $u_{n}$ of the Jacobi  equation  \e{eq:Jy} for $n\to\infty$.  Then we find a representation for the resolvents of the operators $J$ in terms of asymptotic coefficients in formulas for $u_{n}$. Our construction here is independent of Proposition~\ref{Neum}  (or the   von Neumann formulas) and of the results of Sect.~2.2.
  
  \subsection{Setting  the problem}

             To motivate the construction below, let us compare descriptions of self-adjoint extensions of minimal symmetric operators for differential and Jacobi operators. Take, for example, an operator $H_{\rm min}=-d^2/ d x^2$ with domain ${\cal D} (H_{\rm min})= C_{0}^\infty ({\Bbb R}_{+})$ in the space $L^2 ({\Bbb R}_{+})$.  Of course, its adjoint operator $H_{\rm min}^*= H_{\rm max}$ as well as self-adjoint extensions can be described by the von Neumann formulas. It is, however, essentially more convenient  to define these operators in terms of boundary conditions at the point $x=0$.  Indeed, the set ${\cal D} (H_{\rm max})$ coincides with the Sobolev class ${\sf H}^2  ({\Bbb R}_{+})$. Functions $f\in{\sf H}^2  ({\Bbb R}_{+})$  have boundary values $f(0)$ and $f' (0)$.
        The domain  ${\cal D} (\clos H_{\rm min})$ consists of $f\in{\sf H}^2  ({\Bbb R}_{+})$ such that $f(0)=f' (0)=0$, and all self-adjoint extensions of $H_{\rm min}$ are determined by  boundary conditions
          \begin{equation}
f' (0)=h f(0) \;  \mbox{where} \; h\in  {\Bbb R} \q \mbox{or} \q f(0)=0.
   \label{eq:ABX3b}\end{equation}
   
   We are looking for a similar description for self-adjoint extensions of the Jacobi operator $J_{\rm min}$. 
   The starting point of our construction is an asymptotic formula
     \begin{equation}
u_{n} = a_{n}^{-1/2}  \big(s_{+}   e^{i\varphi_{n}} + s_{-}  e^{-i \varphi_{n}}+ o(1)\big)      , \q n\to\infty,
\label{eq:gen}\end{equation} 
for vectors $u\in\{ u_{n}\}\in {\cal D} (J_{\rm max})$.  
 Here the   phases $\varphi_{n} $ are given by an explicit formula (see \e{eq:Grf}, below).    The asymptotic coefficients $s_+= s_+(u)$ and $s_-= s_- (u)$ in \e{eq:gen} depend on $u$.  They play the role of the boundary values $f(0)$ and $f' (0)$ for functions $f\in{\sf H}^2  ({\Bbb R}_{+})$, and we are looking for an analogue of the mapping
     \begin{equation}
f\mapsto(f(0), f' (0))
\label{eq:f}\end{equation}
of the set  ${\cal D} (H_{\rm max})$ onto ${\Bbb C}^2$.

Our  first goal   is to describe all self-adjoint extensions   of  operators  $J_{\rm min}$ in terms of asymptotic coefficients $s_\pm$ in formula \e{eq:gen}. It turns out that such extensions   ${\bf J}_{\omega}$  can be  parametrized by complex numbers $\omega\in{\Bbb T}\subset {\Bbb C}$, and the domain  $ {\cal D} (  {\bf J}_\omega)$ of ${\bf J}_{\omega}$ consists of elements $u\in {\cal D} (    J_{\max})$  satisfying  the condition    
     \begin{equation}
   s_{+} (u) =\omega s_{-} (u) , \q |\omega|=1 .
      \label{eq:gen2}\end{equation}
      This formula plays the role of \e{eq:ABX3b}.

    %  The domain of the operator $\clos J_{\rm min}$ is characterized by the condition $u_{n} = o( \varkappa_{n})$.

   Another goal of this section is to find an explicit formula for the resolvents ${\bf R}_{\omega} (z)=({\bf J}_{\omega}-zI)^{-1}$ of the operators ${\bf J}_{\omega}$. Of course this formula has the same structure as \e{eq:RES1},
      \begin{equation}
{\bf R}_{\omega} (z) h = \boldsymbol{\gamma}_{\omega} (z) \la h,  p(\bar{z})\ra p(z)+ {\cal R}(z)h,
      \label{eq:ABY1}\end{equation}
   but the coefficient
        \begin{equation}
\boldsymbol{\gamma}_{\omega}  (z)= -\frac{s_{+}(q(z))-\omega s_{-}(q(z))} {s_{+}(p(z))-\omega s_{-}(p(z))}
    \label{eq:ABY2}\end{equation}
    is  expressed in terms of the asymptotic coefficients in \e{eq:gen}. 
    
    \subsection{Jost solutions}
    
    We here suppose that the coefficients $a_{n}\to\infty$ as $n\to\infty$  sufficiently rapidly so that condition \e{eq:noCarl} is satisfied (otherwise the operator $J_{\rm min}$ is essentially self-adjoint).  We also assume that  diagonal elements $b_{n}$ are small compared to $a_{n}$. Note that the operator $J_{\rm min}$ is essentially self-adjoint if $b_{n}$ are large compared to $a_{n}$. The critical case where $b_{n}$ and $a_{n}$ are of the same order was considered in \cite{crit-nC}. It is out of the scope of the present paper.

      From analytic point of view, we rely on  the results of paper \cite{nCarl}  where  so called Jost solutions $f_{n}^{(\pm)}(z)$ were distinguished by their asymptotics as $n\to\infty$.  Their construction requires some assumptions on the recurrence coefficients.  In addition to \e{eq:noCarl}, we assume that
      \begin{equation}
   -\frac{b_{n}}{2\sqrt{a_{n-1}a_{n}} }=:\beta_{n}\to  \beta_{\infty}  \; \mbox{where}\;  |\beta_{\infty}| < 1 \;\mbox{as}\;  n\to\infty.
\label{eq:Haupt}\end{equation}
We also require   certain regularity of the behavior of the sequences $a_{n}$ and $\beta_{n}$.    Let us set
   \[
\alpha_{n}=\sqrt{\frac{a_{n+1}}{a_{n }}}, \q k_{n} = \frac{  \alpha_{n-1}} {  \alpha_{n}} =\frac{a_{n}  } {\sqrt{a_{n-1} a_{n+1}} } .
\]
With respect to $a_{n}$, we assume that
\begin{equation}
 \{k_{n} -1 \}\in \ell^1 ({\Bbb Z}_{+} ).
\label{eq:Gr6b}\end{equation} 
With respect to  the ratios $\beta_{n}$, we assume that
     \begin{equation}
  \{\beta_{n}' \}\in \ell^1 ({\Bbb Z}_{+} ).
\label{eq:Gr8}\end{equation}
It is shown in \cite{nCarl}, Lemma~2.1, that
under assumption \e{eq:Gr6b} there exists a finite limit
 \[
\lim_{n\to\infty}  \alpha_{n} = : \alpha_\infty\q \mbox{and} \q\alpha_\infty\geq 1.
\]
 
%  \begin{equation}
%\{\alpha_{n} '  \}\in  \ell^1 ({\Bbb Z}_{+}) .
%\label{eq:Gs9b}\end{equation}

Let us  set  
 \begin{equation}
 \theta_{n}=\arccos\beta_n \q\mbox{and}\q
\varphi_{n}= \sum_{m=0}^{n-1} \theta_{m}, \q n\geq 1,
\label{eq:Grf}\end{equation}
where the sum is restricted to $m$ such that $|\beta_{m}|\leq 1$.  Note that
\[
\varphi_{n}= n\arccos\beta_\infty + o(n)\q \mbox{as}\q n\to\infty.
\]

    \begin{theorem}\label{GSS}\cite[Theorem~4.1]{nCarl}
      Let  assumptions  \e{eq:noCarl},  \e{eq:Haupt}    as well as  \e{eq:Gr6b},  \e{eq:Gr8} be satisfied.
   Then, for all $z\in{\Bbb C} $,   equation \e{eq:Jy} has  solutions $f^{(\pm)}( z ) = (f_{n}^{(\pm)}( z ))$ with asymptotics
    \begin{equation}
f_{n}^{(\pm)}(  z )  = \frac{1}{\sqrt{a_{n}} } e^{\pm i \varphi_{n}} \big(1 + o(1)\big) , \q n\to \infty.
\label{eq:A22G}\end{equation} 
The Wronskian \e{eq:Wr} of these solutions equals
   \begin{equation}
\{ f^{(+)}(  z ), f^{(-)}(  z ) \}=  -2i \alpha_{\infty}^{-1} \sqrt {1-\beta_{\infty}^2}\neq 0,
\label{eq:A2CG1}\end{equation}
so that they    are linearly independent.
For all $n\in {\Bbb Z}_{+}$, the functions $f_{n}^{(\pm)} ( z )$ are entire functions of $z\in  \Bbb C$.
 \end{theorem}
 
 We emphasize that the right-hand side of \e{eq:A22G} depends on $z$ only through the remainder $o(1)$.  Note also that
    \begin{equation}
f_{n}^{(-)}(  z )  = \ov{f_{n}^{(+)}(  \bar{z} )}.
\label{eq:f+-}\end{equation}
 
   \begin{corollary}\label{GSSc}
   Since $f^{(\pm)}(z)\in \ell^2 ({\Bbb Z}_{+})$, under the assumptions of Theorem~\ref{GSS} the operator $J_{\rm min}$ is not essentially self-adjoint. It has deficiency indices $(1,1)$.
 \end{corollary}
 
  \begin{remark}\label{s-adj}
  If the limit in \e{eq:Haupt}  exists but $|\beta_{\infty}|>1$, then,  as shown in \cite{nCarl} (see Theorem~4.11 and Corollary~4.12), the operator $J_{\rm min}$ is   essentially self-adjoint. In the critical case $|\beta_{\infty}|=1$  its essential self-adjointness depends (see \cite{crit-nC}) on details of a behavior of the coefficients $a_{n}$ and $b_{n}$ for $\to\infty$.
 \end{remark}

       \begin{example}\label{GKy}
        Both  conditions  \e{eq:noCarl} and \e{eq:Gr6b}
 are   satisfied for $a_{n}=   n^p$ where   $p>1$ and for $a_{n}=   x^{ n^q}$ where $x>1$, $q< 1$.  In these cases $\alpha_\infty=1$. For $a_{n}=  x^{ n}$,  conditions  \e{eq:noCarl} and \e{eq:Gr6b}
 are   also satisfied but    $\alpha_\infty=\sqrt{x}$. On the  contrary, condition  \e{eq:Gr6b} fails if $a_{n}= x^{ n^q}$  with $q>1$.
      \end{example}
 
  \subsection{Orthogonal polynomials}

An arbitrary solution   of the Jacobi equation  \e{eq:Jy} is a linear combination of the Jost solutions $f_{n}^{(+)}(  z )$ and $f_{n}^{(-)}(  z )$. In particular, this is true for polynomials $ p_{n}(z) $ and $ q_{n}(z) $ so that   
  \begin{equation}
p_{n}(z)= \sigma_{+} (z)f_{n}^{(+)}(  z ) + \sigma_{-} (z) f_{n}^{(-)}(  z )
\label{eq:LC}\end{equation}
and
  \begin{equation}
q_{n}(z)= \tau_{+} (z)f_{n}^{(+)}(  z ) + \tau_{-} (z) f_{n}^{(-)}(  z ),
\label{eq:LCq}\end{equation}
where the coefficients $ \sigma_{\pm} (z)$ and $ \tau_{\pm} (z)$
 can be expressed via the Wronskians:
   \[
  \sigma_{+}(z)= \alpha_{\infty} \frac{ \{p(z), f^{(-)}(  z )\}} {2i \sqrt {1-\beta^2_{\infty}}}, \q
  \sigma_{-}(z)= - \alpha_{\infty} \frac{ \{p(z), f^{(+)}(  z )\}}{2i \sqrt {1-\beta^2_{\infty}}} 
\]
and 
 \[
  \tau_{+}(z)= \alpha_{\infty} \frac{ \{q(z), f^{(-)}(  z ) \} } {2i \sqrt {1-\beta^2_{\infty}}}, \q
  \tau_{-}(z)= - \alpha_{\infty} \frac{\{ q(z), f^{(+)}(  z )\} }{2i \sqrt {1-\beta^2_{\infty}}} .
\]

Observe that
\[
\sigma_{-} (z)= \ov{\sigma_{+}(\bar{z})}\q \mbox{and}\q \tau_{-} (z)= \ov{\tau_{+}(\bar{z})}
\]
 because $p_{n}(z)= \ov{p_{n}(\bar{z})}$, $q_{n}(z)= \ov{q_{n}(\bar{z})}$  and $ f_{n}^{(\pm)}(  z )$ satisfy   \e{eq:f+-}. 
Of course, all coefficients $ \sigma_{\pm}(z)$ and $ \tau_{\pm}(z)$  are  entire functions of $z$.

According to \e{eq:LC}  and \e{eq:LCq}   the following result is a direct consequence of Theorem~\ref{GSS}. 

 \begin{theorem}\label{LC}
    Under the assumptions of Theorem~\ref{GSS} the orthogonal polynomials $p_{n}(z)$ and $q_{n}(z)$
 have   asymptotics, as $n\to\infty$,
  \begin{equation}
p_{n} (z)= a_{n} ^{-1/2} \big(\sigma_{+} (z)e^{ i \varphi_{n}} + \sigma_{-}(z) e^{- i \varphi_{n}}  + o(1)\big)      
\label{eq:LC1P}\end{equation} 
and
 \begin{equation}
q_{n} (z)= a_{n} ^{-1/2} \big(\tau_{+} (z)e^{ i \varphi_{n}} + \tau_{-} (z) e^{-i \varphi_{n}}+ o(1)\big)   .
\label{eq:LC2q}\end{equation} 
   \end{theorem}

    In view of conditions \e{eq:Pz}  and \e{eq:Qz} the Wronskian $\{ p(  z ), q(z)\} =1$. On the other hand, we can calculate this Wronskian using relations   \e{eq:A2CG1} and \e{eq:LC}, \e{eq:LCq}. This yields an identity
   \begin{equation}
  2i \alpha_{\infty}^{-1} \sqrt {1-\beta_{\infty}^2} (\sigma_{+} (z) \tau_{-} (z) -\sigma_{-} (z)\tau_{+} (z))=1,\q \forall z\in{\Bbb C}.
\label{eq:Wro}\end{equation}
We also note an identity  
  \begin{equation}
|\sigma_{+} (z)|^{2} - |\sigma_{-} (z)|^{2}  =\Im z \: \alpha_{\infty}  (1-\beta_{\infty}^2)^{-1/2}  \sum_{n=0}^{\infty}  |p_{n} (z)|^{2}
 \label{eq:LC2p}\end{equation} 
  established in Theorem~4.4 of \cite{nCarl}.

Next, we extend  asymptotic formulas of Theorem~\ref{LC} to all vectors in $ {\cal D}(J_{\rm max})$.  Recall that the number $\Gamma (z; h)$ was defined by formula \e{eq:qres3}.

  \begin{theorem}\label{LCR}
  Let $u= (u_{n})\in{\cal D}(J_{\rm max})$. 
  Choose some $z\in {\Bbb C}$.
     Under the assumptions of Theorem~\ref{GSS} a  sequence $ u_{n} $ has asymptotics
  \begin{equation}
u_{n} = a_{n} ^{-1/2} \big(s_{+}   e^{i \varphi_{n}} +s_{-}  e^{-i \varphi_{n}}+ o(1)\big)      , \q n\to\infty,
\label{eq:LC2}\end{equation} 
where the  coefficients $s_{\pm} = s_{\pm} (u) $ can be constructed by relations
 \begin{equation}
  \begin{split}
s_{+}(u)=  \Gamma(z; ({\cal J}-zI)u) \sigma_{+} (z)+ \la ( {\cal J} -zI) u, p(\bar{z})\ra \tau_{+} (z),
  \\
 s_{-}(u)= \Gamma(z; ({\cal J}-zI)u) \sigma_{-} (z)+ \la ({\cal J}-zI) u, p(\bar{z})\ra \tau_{-} (z).
  \end{split}
  \label{eq:LC31}\end{equation}  

  Conversely, for arbitrary $ s_{+},s_{-}  \in {\Bbb C} $, there exists a vector  
$u \in{\cal D}(J_{\rm max})$ such that    asymptotics \e{eq:LC2}  holds.
 \end{theorem}
 
 %  $\varkappa_{+} (u)=\varkappa_{+}$ and $\varkappa_{-} (u)=\varkappa_{-}$  so that
 
   \begin{pf}
   According to Corollary~\ref{res1}  a vector   $u \in{\cal D}(J_{\rm max})$  admits representation  \e{eq:qres3} where the operator ${\cal R} (z)$  is defined by equality \e{eq:RR11}.
      In view of   relation \e{eq:Ras}  and asymptotics \e{eq:LC1P}, \e{eq:LC2q} we have
    \begin{equation}
 ( {\cal R} (z)h)_{n} =    a_{n} ^{-1/2}  \big( \tau_{+} (z)e^{i \varphi_{n}} + \tau_{-} (z) e^{-i \varphi_{n}} \big) \la h, p(\bar{z})\ra     + o(a_{n} ^{-1/2}), \q n\to\infty,
\label{eq:RR1r}\end{equation}
  for all vectors $h \in \ell^2 ({\Bbb Z}_{+})$.  Therefore it follows from \e{eq:qres3}  that
    \begin{multline*}
u_{n} = a_{n} ^{-1/2}  \Gamma(z; ({\cal J}-zI)u)  \big(\sigma_{+} (z) e^{i \varphi_{n}} + \sigma_{-} (z)e^{-i \varphi_{n}}\big) 
\\
+ a_{n} ^{-1/2} \big(\tau_{+} (z) e^{i \varphi_{n}} + \tau_{-} (z)e^{-i \varphi_{n}}\big) \la  ( {\cal J}-zI) u, p(\bar{z})\ra    + o(a_{n} ^{-1/2})   
 \end{multline*} 
as $ n\to\infty$. This yields  relation \e{eq:LC2} with the coefficients $s_{\pm}$ defined by \e{eq:LC31}.

Conversely, given $ s_{+} $ and $ s_{-} $ and fixing some $z\in{\Bbb C}$, we consider   a system of equations 
 \begin{equation}
  \begin{split}
s_{+}=  \Gamma\sigma_{+} (z)+ \la h , p(\bar{z})\ra \tau_{+} (z),
  \\
 s_{-} = \Gamma \sigma_{-} (z)+ \la h , p(\bar{z})\ra \tau_{-} (z).
  \end{split}
  \label{eq:LC3x}\end{equation}  
   for  $\Gamma$ and $ \la h, p(\bar{z})\ra $. According to \e{eq:Wro} the determinant  of this system is not zero so that  $\Gamma$ and $ \la h, p(\bar{z})\ra $ are uniquely determined by  $s_{+} $ and $ s_{-} $.
 Then we take any $h$ such that its scalar product with $ p(\bar{z})$ equals the found value of  $ \la h, p(\bar{z})\ra$.
  Finally, we define $u$ by formula  \e{eq:qres3}. Asymptotics as $n\to\infty$ of $p_{n}(z)$   and $( {\cal R} (z)h)_{n} $ are given by formulas \e{eq:LC1P} and \e{eq:RR1r}, respectively. In view of equations \e{eq:LC3x} this leads to asymptotics \e{eq:LC2}.
       \end{pf}
       
       %   if and only if, for   all $ z\in {\Bbb C} $,  it admits representation \e{eq:qres3}    where $\Gamma\in{\Bbb C}$ and $h \in \ell^2 ({\Bbb Z}_{+})$ is given by equality  \e{eq:qres2}.  Fix some $z$. 
       
      Theorem~\ref{LCR}  yields a mapping $ {\cal D}(J_{\rm max})\to{\Bbb C}^2$ defined by the formula
       \begin{equation}
u\mapsto (s_{+}  (u) , s_{-}(u)).
\label{eq:mapping}\end{equation}
      The construction of      Theorem~\ref{LCR} depends on the choice of $z\in{\Bbb C}$, but this mapping is defined intrinsically.  In particular, we can set $z=0$ in all formulas of  Theorem~\ref{LCR}.  Note that mapping \e{eq:mapping} is surjective.

Evidently,   \e{eq:mapping} plays the role of mapping \e{eq:f}  for the differential operator $-d^2/ dx^2$ in the space $L^2 ({\Bbb R}_{+})$ and formula \e{eq:qres4} plays the role of the integration-by-parts formula
\[
\int_{0}^\infty f'' (x) \ov{g} (x) dx-\int_{0}^\infty f (x) \ov{g''} (x) dx= f(0) \bar{g}' (0)- f'(0) \bar{g} (0).
\]
 
     Under the assumptions of Theorem~\ref{GSS} the right-hand side of \e{eq:qres4} can be expressed  in terms of the coefficients $s_{+}$ and $s_{-}$.

 \begin{proposition}\label{resAS}
     For all $u, v \in {\cal D} (J_{\rm max})$, we have an identity
    \begin{equation}
\la J_{\rm max}u,v \ra - \la u, J_{\rm max}v \ra = 2 i \alpha_{\infty}^{-1} \sqrt{ 1-\beta_{\infty}^2}\big( s_{+}(u) \ov{s_{+}(v) }- s_{-}(u) \ov{s_{-}(v) }\big)  .
\label{eq:ABS}\end{equation}
 \end{proposition} 
   
    \begin{pf}
    It follows from formula \e{eq:LC2}  that
       \begin{multline*}
\sqrt{a_{n+1}  a_{n } } (u_{n +1}\bar{v}_{n }- u_{n }\bar{v}_{n +1})= \big(s_{+}(u)   e^{ i \varphi_{n+1}} + s_{-} (u)  e^{-i\varphi_{n+1}}\big)\big(\ov{s_{+}(v) }  e^{-i \varphi_{n}} + \ov{ s_{-} (v)}  e^{ i \varphi_{n}}\big)
\\
-  \big(s_{+}(u)   e^{i \varphi_{n}} + s_{-} (u)  e^{-i \varphi_{n}}\big)\big(\ov{s_{+}(v) }  e^{- i \varphi_{n+1}} + \ov{ s_{-} (v)}  e^{ i  \varphi_{n+1}}\big) + o(1)
\\
=\big( s_{+}(u) \ov{s_{+}(v) }-s_{-}(u) \ov{s_{-}(v) }\big) (e^{i\theta_{n}} -e^{-i\theta_{n}}) + o(1).
 \end{multline*}
Passing here to the limit $n\to\infty$ and using equality \e{eq:qres4}, we obtain identity \e{eq:ABS}.
       \end{pf}
       
    We can now characterize set \e{eq:MIN}.
       
      % the domain of the closure 
     %  \[  \clos J_{\min}= J_{\min}^{**}  = J_{\max}^{*}  \]
     %  of the operator $J_{\min}$.
       
       \begin{proposition}\label{ABS}
       A vector
    $v  \in {\cal D} (J_{\rm max})$ belongs to     $ {\sf D}  $ if and only if
    $v_{n}= o(a_{n}^{-1/2})$, that is, 
     \begin{equation}
    s_{+} (v) = s_{-} (v) =0.
      \label{eq:ABS4}\end{equation}
    \end{proposition} 
     
     \begin{pf}    
     A vector $v$ belongs to     $ {\cal D} ( J_{\max}^{*})$ if and only if
       \begin{equation}
    \la J_{\max}u,v \ra= \la u, J_{\max}v \ra
  \label{eq:ABS2}\end{equation}
    for all   $u  \in {\cal D} (J_{\rm max})$.  According to Proposition~\ref{resAS}  equality \e{eq:ABS2} is equivalent to
     \begin{equation}
 s_{+}(u) \ov{s_{+}(v) }- s_{-}(u) \ov{s_{-}(v) }=0.
  \label{eq:ABS3}\end{equation}
  This  is of course true if \e{eq:ABS4} is satisfied. Conversely, if  \e{eq:ABS3} is satisfied for all  $u\in {\cal D} ( J_{\max})$, we use that according to Theorem~\ref{LCR}  the numbers $ s_{+}(u)$ and $s_{-}(u)$  are arbitrary. This implies \e{eq:ABS4}.
         \end{pf}       
    
This result shows that  \e{eq:mapping} considered as a mapping of the factor space  $  {\cal D} (J_{\rm max})/ {\sf D}$ onto ${\Bbb C}^2$ is injective.
  
    \subsection{Self-adjoint extensions}
    
    All self-adjoint extensions ${\bf J}_{\omega} $ of the operator $J_{\min}$ are now parametrized by complex numbers $\omega\in{\Bbb T}\subset {\Bbb C}$.  Let a set $ {\cal D} (  {\bf J}_\omega)\subset {\cal D} (    J_{\max})$ of vectors $u$  be distinguished by  condition     \e{eq:gen2}.

      \begin{theorem}\label{ABX}
           Let the assumptions of Theorem~\ref{GSS}  be satisfied. Then
     for all $\omega\in{\Bbb T}$, the operators ${\bf J}_{\omega}$ are  self-adjoint. Conversely, every   operator $J$ such that 
     \begin{equation}
  J_{\min} \subset  J =J^*\subset J_{\max} 
      \label{eq:ABX2}\end{equation}
      equals ${\bf J}_{\omega}$ for some $\omega\in{\Bbb T}$.
    \end{theorem} 
    
    \begin{pf}
    We proceed from Proposition~\ref{resAS}.   If $u,v \in  {\cal D} (  {\bf J}_\omega)$,  it follows from condition \e{eq:gen2} that   $s_{+}(u) \ov{s_{+}(v) }= s_{-}(u) \ov{s_{-}(v) }$. Therefore   according to equality \e{eq:ABS} $({\bf J}_{\omega}u,v)= (u, {\bf J}_{\omega}v)$ whence ${\bf J}_{\omega}\subset {\bf J}_{\omega}^*$. If $ v\in {\cal D} (  {\bf J}_\omega^*)$, then 
   $\la {\bf J}_{\omega}u,v \ra= \la u, {\bf J}_{\omega}^* v \ra $
    for all $ u\in {\cal D} (  {\bf J}_\omega)$ so that in view of \e{eq:ABS} equality \e{eq:ABS3} is satisfied.
    Therefore $s_{-}(u) (\omega\ov{s_{+}(v) }- \ov{s_{-}(v) })=0$. Since $s_{-}(u)$ is arbitrary, we see  that  $ \omega\ov{s_{+}(v) }- \ov{s_{-}(v) }=0$, and hence $ v\in {\cal D} (  {\bf J}_\omega)$.
    
    Suppose that an operator $J$  satisfies conditions  \e{eq:ABX2}.  Since  $J$ is symmetric, it follows from  Proposition~\ref{resAS} that  equality   \e{eq:ABS3} equality is true  for all $ u, v\in {\cal D} (  J )$.  Setting here $u=v$, we see that $|  s_{+}(v)|=|   s_{-}(v)|$.   There exists a vector $v_{0} \in {\cal D} (  J )$ such that
    $s_{-}(v_{0})\neq 0$ because  $J\neq \clos J_{\rm min}$. Let us set $\omega= s_{+}(v_{0})/ s_{-}(v_{0})$. Then $|\omega|=1$  and relation  \e{eq:gen2}  is a direct consequence of \e{eq:ABS3}.
        \end{pf}

       \subsection{Resolvent}
     
     Now it easy to construct the   resolvent  of the operator ${\bf J}_\omega$  defined in the previous subsection.
     We previously note that, by definition \e{eq:LC2}, 
     \[
     s_{\pm} (p(z))= \sigma_{\pm} (z) \q \mbox{and} \q      s_{\pm} (q(z))= \tau_{\pm} (z).
     \]
     
     \begin{theorem}\label{ABX1}
            Let the assumptions of Theorem~\ref{GSS}  be satisfied. Then
     for all $z\in {\Bbb C}$ with $\Im z\neq 0$ and all $h\in\ell^2 ({\Bbb Z}_{+})$, the resolvent ${\bf R}_{\omega} (z)= ({\bf J}_{\omega}-zI)^{-1}$ of the operator ${\bf J}_\omega$ is given by   equality \e{eq:ABY1} 
      where $\boldsymbol{\gamma}_{\omega}(z)$ is defined by formula  \e{eq:ABY2}, that is,
           \begin{equation}
\boldsymbol{\gamma}_{\omega}(z)= -\frac{\tau_{+}(z)-\omega \tau_{-}(z)} {\sigma_{+}(z)-\omega \sigma_{-}(z)}.
    \label{eq:gam}\end{equation}
    \end{theorem} 
    
         \begin{pf}
           According to      Corollary~\ref{res1} a vector  $u={\bf R}_{\omega}(z) h$ is given by formula \e{eq:qres3} where the coefficient $\Gamma$ is determined by condition 
            \e{eq:gen2}. It follows from Theorem~\ref{LCR}  than the components $u_{n}$ of $u$ have  asymptotics \e{eq:LC2}  with the coefficients $s_{\pm}$ defined by relations 
            \e{eq:LC3x}.  Thus, 
          $u\in {\cal D}({\bf J}_{\omega})$ if and only if  
    \[
    \Gamma \sigma_{+} (z) + \tau_{+} (z) \la h,  p(\bar{z})\ra=\omega \big( \Gamma \sigma_{-} (z) + \tau_{-} (z)\la h,  p(\bar{z})\ra\big)
    \]
    whence
    \[
\Gamma= -\frac{\tau_{+}(z)-\omega \tau_{-}(z)} {\sigma_{+}(z)-\omega \sigma_{-}(z)}\: \la h,  p(\bar{z})\ra.
      \]
Substituting this expression into \e{eq:qres3}, we arrive at  formulas
    \e{eq:ABY1},   \e{eq:ABY2}.
          \end{pf}

          It follows from formula \e{eq:gam} that the spectrum of the operator ${\bf J}_{\omega}$ consists of the points  $z$ where
          \begin{equation}
          \sigma_{+}(z)-\omega \sigma_{-}(z)=0.
         \label{eq:ABY3}\end{equation}
                 Since the functions $\sigma_{+} (z)$ and $\sigma_{-} (z)$  are analytic, the set of such points $z$ is discrete. Moreover,  according to \e{eq:LC2p} $\sigma_{+} (z)\neq \sigma_{-} (z)$ if $\Im z\neq 0$, and therefore zeros $z$ of equation \e{eq:ABY3} lie on the real axis.  This results has  of course to be expected since $z$ are eigenvalues of the self-adjoint operator ${\bf J}_{\omega}$.  We finally note that discreteness of the spectrum of the operators ${\bf J}_{\omega}$ is quite natural because their domains ${\cal D}({\bf J}_{\omega})$ are distinguished by boundary conditions at the point $n=0$ 
 and for $n\to\infty$. Therefore  ${\bf J}_{\omega}$ acquire some features of regular operators.

  \end{document}